\documentclass[12pt]{amsart}
\usepackage{amscd}
\usepackage{xypic}  %commu
\usepackage{amssymb}
\usepackage{amsthm}
\makeindex

\newtheorem{thm}{Theorem}[section]   
\newtheorem{cor}[thm]{Corollary}
\newtheorem{lemma}[thm]{Lemma}
\newtheorem{prop}[thm]{Proposition}
\newtheorem{example}[thm]{Example}
\newtheorem{defn}[thm]{Definition}
\newtheorem{rem}[thm]{Remark}

\newtheorem{lemma-remark}[thm]{Lemma-Remark}
\renewcommand{\proofname}{Proof}

\def\ker{\operatorname{ker}}

\def\min{\operatorname{min}}

\def\im{\operatorname{im}}
\def\rk{\operatorname{rk}}
\def\max{\operatorname{max}}

\def\c1{\operatorname{c_1}}
\def\c2{\operatorname{c_2}}

\def\CC{{\mathbb C}}
\def\ZZ{{\mathbb Z}}

\def\PP{{\mathbb P}}
\def\L{{\mathcal L}}

\def\N{{\mathcal N}}
\def\O{{\mathcal O}}
\def\I{{\mathcal J}}

\def\E{{\mathcal E}}

\def\F{{\mathcal F}}

\def\M{{\mathcal M}}

\def\T{{\mathcal T}}

\def\x{\times}                   % product (fiber)
\def\iso{\simeq}
\def\eqv{\equiv}
\def\sub{\subseteq}

\def\+{\oplus}                   % direct sum
\def\*{\otimes}                  % tensor product
\def\hpil{\longrightarrow}       % ----->

\def\mod{\operatorname{mod}}

\def\Pic{\operatorname{Pic}}

\def\hs{\hspace{.05in}}

\begin{document}

  \title[Curves in families of Calabi-Yau threefolds]{Smooth, isolated curves in families of Calabi-Yau threefolds in homogeneous
         spaces}
  \author{Andreas Leopold Knutsen}

  \address{Department of Mathematics \\
  Johannes Brunsgate 12 \\ 5008 Bergen \\ Norway}
\email{andreas.knutsen@math.uib.no}
\keywords{Isolated curves, deformations, $BN$ general $K3$ surfaces, Calabi-Yau threefolds}
\subjclass{14J32 (14J28)}

\begin{abstract}
 We show the existence of smooth isolated curves of different
 degrees and genera in Calabi-Yau threefolds that are complete
 intersections in homogeneous spaces. Along the way, we classify all degrees and genera of smooth curves on $BN$ general $K3$ surfaces of genus $\mu$, where $5 \leq \mu \leq 10$. By results of Mukai, these are the $K3$ surfaces that can be realised as complete intersections in certain homogeneous spaces. 
\end{abstract}

\maketitle

\section{Introduction}
\label{intro}

In this note we extend the study of embeddings of complex projective curves
into Calabi-Yau complete intersection threefolds from \cite{kn-TAMS} to Calabi-Yau threefolds that are complete intersections in homogeneous spaces. 
For the background and some history of the study of curves in Calabi-Yau threefolds we refer to the introduction of \cite{kn-TAMS}. 

We will pay special attention to families of
Calabi-Yau threefolds in $\PP ^m$, for $7 \leq m \leq 12$, that
are complete intersections in certain homogeneous spaces. We now
briefly fix some notation and refer to Mukai's papers \cite{mu1,mu2} for further details.
We use the following notation: For a vector space $V^i$ of dimension
$i$, we write $G(r,V^i)$ (resp. $G(V^i,r)$) for the Grassmann variety of $r$-dimensional subspaces (resp. quotient spaces) of $V$.
The variety $\Sigma ^{10} _{12} \sub \PP ^{15}$ is a $10$-dimensional spinor variety of degree $12$. Let $V^{10}$ be a $10$-dimensional vector space with a nondegenerate second symmetric tensor $\lambda$. Then $\Sigma ^{10} _{12}$ is one of the two components of the subset of $G(V^{10}, 5)$ consisting of $5$-dimensional totally isotropic quotient spaces. 
The variety $\Sigma ^6 _{16} \sub \PP ^{13}$ is the Grassmann variety of $3$-dimensional totally isotropic quotient spaces of a $6$-dimensional vector space $V^6$ with a nondegenerate second skew-symmetric tensor $\sigma$. It has dimension $6$ and degree $16$.
Finally, $\Sigma ^5 _{18}= G/{\mathcal P} \sub \PP ^{13}$, where $G$ is the automorphism group of the Cayley algebra over $\CC$ and $\mathcal{P}$ is a maximal parabolic subgroup. The variety has dimension $5$ and degree $18$.

By adjunction, we get the following families of Calabi-Yau threefolds
in $\PP^m$
as complete intersections in homogeneous spaces (we use the notation
$(a_1,...,a_n) \cap \Sigma$ for a complete intersection of type
$(a_1,...,a_n)$ in $\Sigma$):

\vspace{.4cm}
\begin{center}
\begin{tabular}{|c|l|} \hline
$m$ & intersection type of Calabi-Yau threefold \\ \hline
$7$  &  \hspace{1.4cm}  $(3,1,1)  \cap G(2, V^5)$     \\    \hline
$8$  &  \hspace{1.4cm}  $(2,2,1)  \cap G(2, V^5)$     \\    \hline
$9$  &  \hspace{1.4cm} $(2,1 ^6) \cap \Sigma ^{10} _{12}$     \\    \hline
$10$  &  \hspace{1.4cm} $(2,1^4)  \cap G(V^6, 2)$     \\    \hline
$11$  &  \hspace{1.4cm} $(2,1,1)  \cap \Sigma ^6 _{16}$     \\    \hline
$12$  & \hspace{1.4cm}  $(2,1)    \cap \Sigma ^5 _{18}$    \\    \hline
\end{tabular}
\end{center}
\vspace{.4cm}

We will prove the existence of isolated smooth curves of various
degrees and (low) genera in each of these Calabi-Yau threefolds. More precisely, we will prove the following result. 

\begin{thm} \label{result}
The general Calabi-Yau threefold of type $(3,1,1) \cap G(2, V^5) \sub \PP^7$
contains an isolated smooth irreducible curve of degree $d$ and  genus $g$ if $0 \leq g \leq 4$ and $d^2 \geq 20g-4$ or  $5 \leq g \leq 13$ and $d \geq g+6$.

The general Calabi-Yau threefold of any of the other types in the table above 
contains an isolated smooth irreducible curve of degree $d$ and  genus $g$ if $0 \leq g \leq m-4$ and $d^2 \geq 4(m-3)g-4$.
\end{thm}

The main purpose of this note is to show how to apply the results from \cite{kn-TAMS} to the case where the ambient variety is not a projective space. In fact, the main result in \cite{kn-TAMS} holds for an arbitrary ambient space, but the main application therein was to the case where the ambient space is a projective space, yielding an existence result similar to Theorem \ref{result} for the five different types of complete intersection Calabi-Yau threefolds. 

We will use the same technique from \cite{kn-TAMS} (and from references therein),
namely we
start with a $K3$  surface of a suitable  intersection type in one of
the homogeneous spaces above containing a particular smooth curve and
construct a (nodal) Calabi-Yau threefold of a corresponding
intersection type in the same ambient space containing the surface. Then
we use the main result from \cite{kn-TAMS} to show that finitely many of these curves deform to the general such Calabi-Yau threefold as smooth and isolated curves.

The main difficulty (and the main reason why we did not include this in \cite{kn-TAMS}) is to know which curves may exist on the $K3$ surfaces in question. For this we use the precise characterization of such surfaces due to Mukai \cite{mu1,mu2} and we obtain a complete answer in Theorem \ref{bncurves}, which is of independent interest and is the second main result of this note. 

Section \ref{sec:proofthm} is devoted to showing that the proof in \cite{kn-TAMS} in the case of a complete intersection in a projective space goes through almost unaffected in this more general setting.

\subsection{Conventions and definitions}
The ground field is the field of complex numbers. We say a curve $C$
in a variety $V$ is {\it infinitesimally rigid} or {\it isolated} in $V$
if the space of
embedded deformations of $C$ in $V$ is zero-dimensional and reduced, equivalently, if $H^0(C, \N _{C/V})=0$.

By a $K3$ surface is meant a reduced and irreducible smooth surface $X$ with trivial
canonical bundle and such that $H^1 (\O_X)=0$. In particular $h^2 (\O_X)= 1$
and $\chi (\O_X)= 2$.

A Calabi-Yau threefold $Y$ is a projective variety of dimension $3$ with
trivial canonical bundle and $h^1(\O_Y)=h^2(\O_Y)=0$. In this paper a
Calabi-Yau threefold will be at worst nodal.

We will use line bundles and divisors on a surface will little or no distinction. In particular, we use the following convention for a divisor (class) or line bundle $D$: We write $D \geq 0$ to mean that $h^0(D)>0$ and $D>0$ if moreover $D$ is nontrivial.

We say that a variety $X$ is a complete intersection in a variety $P \subset \PP^m$  if there are hypersurfaces $H_1,\ldots,H_k$ in $\PP^m$ such that
$X=P \cap H_1 \cap \cdots \cap H_k$ scheme-theoretically.

\subsection{Notice} This note is a revised version of a first version with approximately the same title from 2001. That version could not be published, since - as was pointed out by a careful referee - it 
used a published result  that turned out to have a serious gap. The gap was only fixed recently under some slightly stronger conditions in \cite{kn-TAMS}, and we refer to that paper for more details. 

The stronger conditions implied that the first version needed a rewrite (in fact, the new version of Theorem \ref{result} is valid in a narrower range of values of $d$ and $g$ than the first version). We also remark that the exposition of Section \ref{BNK3} of this revised version is considerably simplified and improved. 

\subsection{Acknowledgements}
The author thanks T.~Johnsen, H.~Clemens
and H.~P.~\linebreak Kley for useful discussions, H.~Murakami for patiently translating relevant passages of \cite{mu2}
from Japanese to English to the author, and S.~Mukai for
providing the author with a preliminary English translation. 
The author also thanks the referee of the 2001 version of this paper for having discovered the gap in the result that was used and for explaining this very carefully  in the report.

A part of this work was done during the author's stay at the
Mittag-Leffler Institute in Djursholm during the academic year 1998-99 and during the author's stay at the Department of Mathematics
at the University of Utah, Salt Lake City, in the Spring of 2000. The author wishes to thank
both institutions for their hospitality and warm atmosphere.

\section{$BN$ general $K3$ surfaces}
\label{BNK3}

A {\it polarized} $K3$ surface of genus $\mu$ (or degree $2(\mu-1)$) is  a pair $(X,H)$ where $X$ is
a $K3$ surface and $H$ is an ample line bundle on $X$ with $p_a(H)=\mu$ (equivalently $H^2=2(\mu-1)$. If in addition $H$ is indivisible in $\Pic X$, then the surface is called {\it primitively polarized}.

\begin{defn}  \label{BNgendef} 
 {\rm (\cite[Def.~3.8]{mu2}) A polarized $K3$ surface $(X,H)$ of genus ${\mu} \geq 2$ is said to be {\rm Brill-Noether ($BN$) general} if the inequality $h^0 (M) h^0(N) < h^0(H)={\mu}+1$ holds for any pair $(M,N)$ of non-trivial line bundles such that $M +N \sim H$.}
\end{defn}

\begin{rem}  \label{BNgensuff} 
 {\rm It is easy to see that the condition of being $BN$ general is satisfied 
if there is a smooth Brill-Noether general curve in $|H|$. (Indeed, if there is an $M$ such that $h^0 (M) h^0(H-M) \geq {\mu}+1$, then the Brill-Noether number $\rho(\O_C(M)) <0$ on any smooth 
$C \in |H|$.) We do not know if the converse is true. Also see Remark \ref{rem:KK3} below.}
\end{rem}

   Since the generic primitively polarized  $K3$ surface $(X,H)$ of genus ${\mu}$
   satisfies $\Pic X \iso \ZZ[H]$, the primitively polarized $K3$ surfaces that are {\it not} $BN$ general form a countable union of proper closed subsets  
in the moduli space of primitively polarized $K3$ surfaces of a fixed genus ${\mu}$.

   \begin{rem}
{\rm By the classical results of Saint-Donat \cite[Thm.~3.1 and Prop.~8.1]{S-D}, ample line bundles on $K3$ surfaces that are not globally generated (resp. very ample) are completely characterised by their self-intersections and by the existence of divisors on $S$ with specific intersection properties. In particular, it is easily checked that if $(S,H)$ is $BN$ general, then $H$ is globally generated, whence it defines a morphism $\varphi_{H}$ into $\PP^{\mu}$, as $h^0(H)=\mu+1$ by Riemann-Roch and Kodaira vanishing. Moreover, if $\mu \geq 3$, one can also check that $H$ is very ample, so that the morphism is an embedding.}
   \end{rem}

The following theorem is due to Mukai and gives a description of $BN$ general $K3$ surfaces of low genus as complete intersections in the homogeneous spaces mentioned in the introduction:

\begin{thm} \label{BNgenthm}
{\rm (\cite{mu1,mu2})} 
The projective models of $BN$ general polarized $K3$
  surfaces of small genus are as in the table below. Conversely, any surface in the table is a $BN$ general $K3$ surface.
\end{thm}

\begin{center}
\begin{tabular}{|c|l|} \hline
genus & projective model of $BN$ general polarized $K3$ surface \\ \hline
$2$  &  \hspace{.5cm} $X_2 \hpil \PP ^2$ double covering with branch sextic  \\  \hline
$3$  &  \hspace{.5cm} $(4) \sub \PP ^3$                \\    \hline
$4$  &  \hspace{.5cm} $(2,3) \sub \PP ^4$              \\    \hline
$5$  &  \hspace{.5cm} $(2,2,2) \sub \PP ^5$            \\    \hline
$6$  &  \hspace{.5cm} $(1,1,1,2) \cap G(2, V^5) \sub \PP ^6$      \\  \hline
$7$  &  \hspace{.5cm} $(1 ^8)\cap \Sigma ^{10} _{12} \sub \PP ^7$ \\  \hline
$8$  &  \hspace{.5cm} $(1^6) \cap G(V^6, 2) \sub \PP ^8$          \\  \hline
$9$  &  \hspace{.5cm} $(1^4) \cap \Sigma ^6 _{16} \sub \PP ^9$    \\  \hline
$10$  & \hspace{.5cm} $(1^3) \cap \Sigma ^5 _{18} \sub \PP ^{10}$ \\  \hline
%$12$  & \hspace{.5cm} $S_{12} =  (1)     \sub \Sigma _{12}$   \\    \hline
\end{tabular}
\end{center}
\vspace{0.6cm}

It is readily checked and very well-known that any polarized $K3$ surface of genus $2$ is automatically $BN$ general, and that a polarized $K3$ of genus $3$ and $4$ is $BN$ general if and only if $H$ is very ample. In particular, any smooth $K3$ surface embedded as a linearly normal surface in $\PP^3$ or $\PP^4$ is a quartic or  a complete intersection of a hyperquadric and a hypercubic respectively. By contrast, if $\mu=5$, not all embedded $K3$s are complete intersections of three hyperquadrics: there is a codimension one subspace of polarized $K3$ surfaces of genus $5$ that are  a section of $|\O_Y(3)-\F|$ in a smooth three-dimensional rational normal scroll $Y$ in $\PP^5$, where $\F$ denotes the class of the $\PP^2$-fibers \cite{S-D,JK}.

In the next section we will study the possible combinations of pairs $(d,g)$ for a smooth curve $C$ of degree $d$ and genus $g$ on a $BN$ general $K3$ surface of genus $\mu$, with $6 \leq \mu \leq 10$. (The case $\mu=5$ has already been covered in \cite[Thm.~6.1(3)]{kn}). 

We will make use of the following result, which in a reformulated way says that in Definition \ref{BNgendef} we may restrict to considering line bundles $M$ and $N$ that are effective, nef and with vanishing $H^1$.

\begin{prop} \label{prop:BN}
A polarized $K3$ surface $(X,H)$ of genus $\mu$ is $BN$ general if and only if for any pair $(M,N)$ of effective, nef line bundles such that $M+N \sim H$, the inequality
\begin{equation} \label{eq:BN}
 \left(\frac{1}{2}M^2+2\right)\left(\frac{1}{2}N^2+2\right) \leq \frac{1}{2}H^2+1=\mu
\end{equation}
holds.
\end{prop}

\begin{proof}
 If $H \sim M+N$ is an effective decomposition, then by Riemann-Roch $h^0(M) \geq \frac{1}{2}M^2+2$ and similarly for $N$. Hence, if \eqref{eq:BN} is not satisfied, we have
\[ h^0(M)h^0(N) \geq \left(\frac{1}{2}M^2+2\right)\left(\frac{1}{2}N^2+2\right)\geq \frac{1}{2}H^2+2=h^0(H), \]
the latter equality following from Riemann-Roch and ampleness of $H$. Hence $(X,H)$ is not $BN$ general.

Conversely, assume that $(X,H)$ is not $BN$ general and let $H \sim M+N$ be  an effective decomposition such that $h^0(M)h^0(N) \geq h^0(H) = \frac{1}{2}H^2+2$. 

We first claim that we can assume that $M$ is nef.

Indeed, let $\Gamma$ be a smooth irreducible curve satisfying $\Gamma.M <0$ (and, necessarily, $\Gamma^2=-2$.) Then $h^0(M-\Gamma)=h^0(M)$ and $h^0(N+\Gamma) \geq h^0(N)$. Hence $h^0(M-\Gamma)h^0(N+\Gamma) \geq h^0(H)$. Repeating the process, and since all $(-2)$-curves we ``remove'' from $M$ are part of the base divisor of $|M|$ (whence are finitely many), we may assume that $M$ is nef. 

We next claim that we can also assume that $N$ is nef. 

Indeed, let  $\Gamma$ be a smooth irreducible curve satisfying $\Gamma.N <0$ and $\Gamma^2=-2$. As $H$ is ample, we must have $\Gamma.M \geq 2$.
Hence $M+\Gamma$ is still nef,  $h^0(N-\Gamma)=h^0(N)$ and $h^0(M+\Gamma) \geq h^0(M)$. Repeating the process, we may assume that $N$ is nef.

In particular, we may assume that $M^2 \geq 0$ and $N^2 \geq 0$ and that both $M$ and $N$ are nontrivial. In particular, $h^2(M)=h^2(N)=0$, and we will be done by Riemann-Roch by proving that we may choose $M$ and $N$ such as $h^1(M)=h^1(N)=0$.

Assume that $h^1(M) \neq 0$. Then, as $M$ is nef, we must have $M \sim lE$ for some integer $l \geq 2$, where $|E|$ is an elliptic pencil \cite{S-D}. 
Pick $l$ distinct smooth members of $|E|$.
From 
\[
\xymatrix{
0  \ar[r] & N \ar[r] & H \ar[r] & H_{|E_1 \sqcup \cdots \sqcup E_l} \ar[r] & 0
}
\]
we find that
\[
h^0(H) = \chi(N)+lE.H = \chi(N)+lE.(lE+N)=\chi(N)+lE.N.
\]
Now $E.N=E.H>0$ by ampleness of $H$, whence $N+(l-1)E$ is big and nef, in particular $h^1(N+(l-1)E)=0$. From
\[
\xymatrix{
0  \ar[r] & N \ar[r] & N+(l-1)E \ar[r] & (N+(l-1)E)_{|E_1 \sqcup \cdots \sqcup E_{l-1}} \cong 
N_{|E_1 \sqcup \cdots \sqcup E_{l-1}}
\ar[r] & 0
}\]
we find that
\[
h^0(N+(l-1)E) = \chi(N)+(l-1)E.N.
\]
Hence
\[h^0(E)h^0(N+(l-1)E)=2(\chi(N)+(l-1)E.N) > \chi(N)+lE.N=h^0(H),\]
and the decomposition $H\sim E+(N+(l-1)E)$ has the desired property.
\end{proof} 

\begin{example}
{\rm Let $\mu=5$ and assume that $H \sim M+N$ with $M,N$ effective, nontrivial and nef. Assume, without loss of generality, that $M^2 \leq N^2$.
It is readily verified, using the Hodge index theorem, that the only possibilities are
$(M^2,N^2,M.N)=(0,0,4)$, $(0,2,3)$, $(0,4,2)$, $(0,6,1)$ and $(2,2,2)$, the latter implying that $M\sim N$.
Inequality \eqref{eq:BN} is satisfied only in the first case. Moreover, the three middle cases can equally be characterized by $M^2=0$ and $M.L \leq 3$. Hence we retrieve Saint-Donat's well-known result \cite[(7.1) and Thm.~7.2]{S-D} that $X$ can be embedded as a complete intersection of three hyperquadrics if and only if there is no effective divisor $D$ on $X$ satisfying $D^2=0$ and $D.H \leq 3$ or $D^2=2$ and $H \sim 2D$. 

The cases where a divisor $D$ as above exist can be described as follows. The cases where there is a divisor satisfying $D^2=0$ and $D.H=1$ are precisely the cases where $H$ is not globally generated \cite[Thm.~3.1]{S-D}. The cases where there is a divisor satisfying $(D^2,D.H)=(0,1)$ or $D^2=2$ and $H \sim 2D$ are 
precisely the cases where $H$ is not very ample \cite[Prop.~8.1]{S-D}. Assuming $H$ is globally generated, the morphism $\varphi_{H}$ maps $X$ generically $2:1$ onto a rational normal scroll in the first case, and onto the Veronese surface in $\PP^5$ in the second case. Finally, assuming $H$ is very ample, the cases where there is a divisor satisfying $D^2=0$ and $D.H=3$ are precisely the cases mentioned above where $\varphi_{H}(X)$ is a section of a smooth three-dimensional rational normal scroll $Y$ in $\PP^5$.

In the same way we find that a polarized $K3$ surface $(X,H)$ of genus $\mu=6$ can be realized as $(1,1,1,2) \cap G(2, V^5) \sub \PP ^6$ if and only if there is no effective divisor $D$ on $X$ satisfying $D^2=0$ and $D.H \leq 3$ or $D^2=2$ and $D.H=5$. Still by \cite[(7.1) and Thm.~7.2]{S-D} this is equivalent to the homogeneous ideal of $X$ being generated only by hyperquadrics (case (ii) in 
\cite[Thm.~7.2]{S-D} is equivalent to the case $D^2=2$ and $D.H=5$.) 

In a similar way we can obtain concrete descriptions of $BN$ general $K3$ surfaces of genus $\mu=7,8,9,10$ in terms of nonexistence of divisors on $X$ with low self-intersections, but we leave this to the reader. For $\mu \geq 7$, however, the characterization is no longer equivalent to Saint-Donat's characterization of the homogeneous ideal of $X$ being generated only by hyperquadrics.
All projective models of $K3$ surfaces of genus $\mu$ with $5 \leq \mu\leq 10$ that are not $BN$ general are given in \cite{JK} as sections in scrolls.
}
\end{example}

\section{Smooth curves on $BN$ general $K3$ surfaces}
\label{smcur}

We first recall the main result in \cite{kn}:

\begin{thm} \label{mainthm}
  {\rm (\cite[Thm. 1.1]{kn})} Let $n \geq 2$, $d>0$, $g \geq 0$ be integers. Then there exists a
  $K3$ surface $X$ of degree $2n$ in $\PP^{n+1}$ containing a smooth irreducible
  curve $C$  of degree $d$ and genus $g$ if and only if
    \begin{itemize}
     \item[(i)] $g=d^2/4n + 1$ and there exist integers $k,m \geq 1$ and $(k,m)
      \not = (2,1)$ such that $n = k^2 m$ and $2n$ divides $kd$,
     \item[(ii)] $d^2/4n < g < d^2/4n +1$ except in the following cases
       \begin{itemize}
        \item[(a)] $d \eqv \pm 1, \pm 2 \hs (\mod 2n)$,
        \item[(b)] $d^2-4n(g-1) = 1$ and $d \eqv n \pm 1 \hs (\mod 2n)$,
        \item[(c)] $d^2-4n(g-1) = n$ and $d \eqv n \hs (\mod 2n)$,
        \item[(d)] $d^2-4n(g-1) = 1$ and $d-1$ or $d+1$ divides $2n$,
       \end{itemize}
     \item[(iii)] $g=d^2/4n$ and $d$ is not divisible by $2n$,
     \item[(iv)] $g < d^2/4n$ and $(d,g) \not =(2n+1,n+1)$.
    \end{itemize}

  Furthermore, in case (i) $X$ can be chosen such that $\Pic X = \ZZ
  \frac{2n}{dk}[C] = \ZZ \frac{1}{k}[H]$ and in cases (ii)-(iv) such that $\Pic
  X = \ZZ [H] \+ \ZZ [C]$, where $H$ is the hyperplane section of $X$.
\end{thm}

As remarked above, the surfaces with $n \leq 3$ are automatically $BN$ general. 

For all possible triples of integers $(n,d,g)$, with $4 \leq 5 \leq 9$, $d >0$ and $g\geq 0$, we will now investigate which triples are also possible for smooth curves on $BN$ general $K3$ surfaces (retrieving the case $n=4$ from \cite[Thm. 6.1(3)]{kn}). The rest of this section is devoted to the proof of the following theorem of independent interest:

\begin{thm} \label{bncurves}
 Let $4 \leq n \leq 9$, $d>0$, $g \geq 0$ be integers. Then there exists a
  $BN$ general $K3$ surface $X$ of degree $2n$ in $\PP^{n+1}$ (and genus $\mu=n+1$)
 containing a smooth irreducible
  curve $C$  of degree $d$ and genus $g$ if and only if we are in one of the following cases:
    \begin{itemize}
     \item[(I)] $g=d^2/4n + 1$ and  $2n$ divides $d$;
     \item[(II)] $n=5$ and $g=(d^2+4)/20$;
\item[(III)] $n=7$ and $g=(d^2+3)/28$;
\item[(IV)] $g=d^2/4n$ and $d$ is not divisible by $2n$;
     \item[(V)] $g < d^2/4n$ and $(d,g) \not =(2n+1,n+1)$.
    \end{itemize}

  Furthermore, in case (I) the surface  $X$ can be chosen such that $\Pic X = \ZZ[H]$ and in cases (II)-(V) such that $\Pic
  X = \ZZ [H] \+ \ZZ [C]$, where $H$ is the hyperplane section of $X$.
\end{thm}

We start with an easy result for the cases $d^2-4n(g-1)=0$.

\begin{lemma} \label{primocaso}
Let $n \geq 2$, $d >0$ and $g >0$ be integers satisfying $d^2-4n(g-1)=0$. Then
there exists a $BN$ general $K3$ surface $X$ of degree $2n$ in $\PP^{n+1}$ 
containing a smooth curve $C$  of degree $d$ and genus $g$ if and only if $2n$
divides $d$. 
\end{lemma}

\begin{proof}
If $2n$ divides $d$, then by Theorem \ref{mainthm} there is a $K3$ surface $X$ of 
degree $2n$ in $\PP^{n+1}$  containing a smooth curve $C$  of
degree $d$ and genus $g$ such that $\Pic X \iso \ZZ[H]$, where $H$ is the
hyperplane section (and $C \sim \frac {d}{2n}H$). Such a surface is clearly
$BN$ general.

Conversely, if $X$ is a $K3$ surface of 
degree $2n$ in $\PP^{n+1}$  containing a smooth curve $C$  of
degree $d$ and genus $g$, then by the Hodge index theorem $2nC \sim dH$, where
$H$ is the hyperplane section. If $2n$ does not divide $d$, then $H$ must be 
divisible, whence it is easily seen that $X$ is not $BN$ general.
\end{proof}

The surfaces in Theorem \ref{mainthm} with $d^2-4n(g-1)>0$ can all be
chosen with
the property that $\Pic X = \ZZ [H] \+ \ZZ [C]$. It is easily seen that we must have $\rk (\Pic X) \geq 2$ in these cases. 

\begin{defn}
 {\rm We say that a triple $(X,H,C)$ is} an $(n,d,g)$-surface {\rm if 
$X$ is a smooth $K3$ surface of degree $2n$ in $\PP^{n+1}$, $H$ its hyperplane bundle, and $C$ an effective 
divisor with $C^2=2g-2$, $C.H=d$  and $d^2-4n(g-1)>0$. (Often we will denote it only by $X$, when $H$ and $C$ are understood.) 

We say that $(X,H,C)$ is} Picard minimal {\rm if furthermore $\Pic X = \ZZ[H] \+ \ZZ[C]$}.
\end{defn}

Even if we are ultimately interested in the cases where $C$ is a smooth and irreducible curve, we do not require that $C$ be smooth and irreducible in the definition. The background for this is a reduction procedure in Lemmas \ref{small} and \ref{lemma:caso2} below where we will substitute $C$ by another effective divisor that may not be smooth and irreducible even if $C$ is. 

\begin{rem} \label{rem:KK3} 
{\rm The Picard minimal $(n,d,g)$-surfaces in our definition with smooth $C$ are precisely the surfaces denoted by $S_{n,d,g}$ and called Knutsen $K3$ surfaces in \cite{AM} 
(for $d^2-4n(g-1)>0$). The main result of \cite{AM} gives necessary and sufficient numerical conditions for such surfaces to have the strong property that all their hyperplane sections   
are reduced and irreducible. In particular, such surfaces are obviously $BN$ general, and furthermore, by Lazarsfeld's famous result \cite{la}, all their hyperplane sections are Brill-Noether general. The interest in this is related to the question whether the surface embeds in a certain Fano threefold, cf. \cite{ACM}.}
\end{rem}

Assume that we have an effective decomposition $H \sim D+F$ for $D>0$ and $F >0$ on a Picard minimal $(n,d,g)$-surface. Then we can write $D = aH-bC$
and $F = (1-a)H + bC$, for some integers $a$, $b$. Since $D$ and $F$
are effective, one easily sees that we can assume that
\[ a \geq 1  \; \mbox{  and  } \; b \geq 1. \]

\begin{lemma} \label{triangle}
 Let $X$ be a Picard minimal $(n,d,g)$-surface.

 If $d \geq n+g$, then there is no nontrivial effective decomposition
  $H \sim D+F$ with $D$ and $F$ both nef. 
In particular, $X$ is $BN$ general.
\end{lemma}

\begin{proof}
Assume that $H \sim D+F$ is an effective decomposition, with $D$ and $F$ nef, and let $a,b$ be as above. 

If $b \geq a$, then $D=aH-bC=b(H-C)-(b-a)H$, so that $b(H-C)$ is effective. 
Since $(H-C)^2=2(n-d+g-1) \leq -2$ by assumption, there is an irreducible curve $\Gamma$ satisfying $\Gamma^2=-2$ and $\Gamma.(H-C) <0$. Hence also $\Gamma.D<0$, contradicting the nefness of $D$. 

If $b<a$ a similar argument contradicts nefness of $F$. 

The last statement follows from Proposition \ref{prop:BN}. 
\end{proof}

\begin{lemma} \label{rational}
 (i) Any Picard minimal $(n,d,0)$-surface is $BN$ general.

 (ii) If $d > \frac{n+2}{2}$, any Picard minimal $(n,d,1)$-surface is $BN$ general. 
     If $d \leq \frac{n+2}{2}$, no $(n,d,1)$-surface is $BN$ general. 
\end{lemma}

\begin{proof}
Let $X$ be any Picard minimal $(n,d,g)$-surface with $g \leq 1$. Assume that $H \sim D+F$ is an effective decomposition, with $D$ and $F$ nef, and let $a,b$ be as above. As $C.F=C.((1-a)H+bC) =(1-a)d + 2b(g-1) \leq (1-a)d$, nefness of $F$ implies that $a=1$ and $g=1$. Hence (i) is proved, by Proposition \ref{prop:BN}.

Let now $g=1$. We have
\begin{equation} \label{eq:ine}
 \left(\frac{1}{2}D^2+2\right)\left(\frac{1}{2}F^2+2\right) =
\left(\frac{1}{2}(H-bC)^2+2\right)\left(\frac{1}{2}(bC)^2+2\right)=2(n-bd+2). 
\end{equation}
If $d > \frac{n+2}{2}$, then $2(n-bd+2) \leq 2(n-d+2) \leq n+1 =\frac{1}{2}H^2+1$, showing that $X$ is $BN$ general by Proposition \ref{prop:BN}.

Finally, assume that $X$ is any $(n,d,1)$-surface and $d \leq \frac{n+2}{2}$.
Let $D=H-C$ and $F=C$. Then 
$2(n-d+2) \geq n+2 = \frac{1}{2}H^2+1$, so that \eqref{eq:ine} (with $b=1$) and  Proposition \ref{prop:BN} show that $X$ is not $BN$ general.
\end{proof}

The next lemma will be useful for considering all remaining cases. The idea is to reduce to smaller values of $d$ and $g$.

\begin{lemma} \label{small}
  Let $n \geq 2$, $d \geq 1$ and $ g \geq 0$ be integers satisfying $g <
  d^2 /4n +1$ and let $(X,H,C)$ be an $(n,d,g)$-surface. Define
\[ k:= \max \left\{i \in \ZZ_{\geq 0} \; | \; d-2ni>0 \; \; \mbox{and} \; \; g-di+ni^2>0 \right\} \]
and set $C_0:=C-kH$, $d_0:=d-2nk$ and $g_0:=g-dk+nk^2$.

Then $C_0>0$, $C_0.H=d_0>0$, $C_0^2=2(g_0-1)\geq 0$ and $d_0^2-4ng_0=d^2-4ng$. Moreover, precisely one of the following two cases occurs:
\begin{itemize}
\item[(i)] $d_0 \geq g_0+n$ and $g_0 \leq \frac{d_0^2}{4n}$;
\item[(ii)] $d_0<2n$ and $d_0-n <g_0<\frac{d_0^2}{4n}+1$.
\end{itemize}
\end{lemma}

\begin{proof}
The numerical properties of $C_0$, $d_0$ and $g_0$ are immediate. The fact that
$C_0 > 0$ follows from Riemann-Roch using the facts that $C_0^2 \geq 0$ and $h^2(C_0)=h^0(-C_0)=0$ as $C_0.H>0$.

The two cases (i) and (ii) are clearly mutually exclusive, so we have left to prove that one of the cases  occurs. 

Assume that we are not in case (i). It is easily verified that the second numerical condition in (i) is a consequence of the first, whence we must have $d_0 <g_0+n$. Then $d_0-n <g_0<\frac{d_0^2}{4n}+1$ is clear, so we have left to prove that $d_0<2n$. Assume, to get a contradiction, that 
 $d_0 \geq 2n$. Then $(C-(k+1)H).H=d_0-2n \geq 0$, whence $h^2(C-(k+1)H)=h^0(-(C-(k+1)H))=0$. Moreover, $(C-(k+1)H)^2=2[g-d(k+1)+n(k+1)^2-1]=2[g_0-d_0+n-1]\geq0$,
so that
$C-(k+1)H >0$ by Riemann-Roch and the fact that $H$ and $C$ are independent in $\Pic X$. Hence $d-2n(k+1)= d_0-2n = (C-(k+1)H).H>0$, contradicting the definition of $k$. 
\end{proof}

Clearly the $(n,d,g)$-surface $(X,H,C)$ is $BN$ general if and only if the 
$(n,d_0,g_0)$-surface $(X,H,C_0)$ is $BN$ general.

If we end up in case (i), we are done by the following:

\begin{cor} \label{cor:small}
  Let $X$ be a Picard minimal $(n,d,g)$-surface with $g <
  d^2 /4n +1$ and $d_0$, $g_0$ be as in Lemma \ref{small}.

If $d_0 \geq g_0+n$, then $X$ is $BN$ general.
\end{cor}

\begin{proof}
Apply Lemma \ref{triangle} to the Picard minimal $(n,d_0,g_0)$-surface $(X,H,C_0)$, where
$C_0$ is as in  Lemma \ref{small}.
\end{proof}

If we end up in case (ii), the following result will be useful:

\begin{lemma} \label{lemma:caso2}
  With the same assumptions and notation as in Lemma \ref{small}, if case (ii) occurs, set $C':=H-C_0$, $d':=2n-d_0$ and $g':=n-d_0+g_0$. Then ${d'}^2-4ng'=d^2-4ng$, $C'>0$, $C'.H=d'>0$ and ${C'}^2=2(g'-1)\geq 0$.
\end{lemma}

\begin{proof}
  Easy computation. 
\end{proof}

We will also need the following:

\begin{lemma} \label{lemma:spec}
  Any Picard minimal $(n,d,g)$-surface with $(n,d,g)=(8,8,2)$ or $(9,9,2)$ is $BN$ general.
\end{lemma}

\begin{proof}
  Assume that $H \sim D+F$ is an effective decomposition, with $D$ and $F$ nef. As usual, write $D := aH-bC$ and $F := (1-a)H + bC$, for two integers
$a$, $b \geq 1$, and assume $D,F >0$. Since $H$ is (very) ample, we have
\[ D.H = 2na-bd >0 \]
and
\[ F.H =(1-a)2n +bd >0, \]
which gives
\begin{equation} \label{eq:bn1}
 \frac{2n}{d}(a-1) < b < \frac{2n}{d}a.
\end{equation}
Since $n=d$, this yields $b=2a-1$. Furthermore, nefness of $F$ yields
\[ 0 \leq F.C =d(1-a)+ 2(g-1)b = d(1-a)+2b=d(1-a)+2(2a-1)= d-2-(d-4)a, \]
which gives the only possibility $a=1$ for $d=8,9$. Hence also $b=1$ and we 
have $D=H-C$ and $F=C$. 

We have 
\begin{equation}
\left(\frac{1}{2}(H-C)^2+2\right)\left(\frac{1}{2}C^2+2\right)  =   (n-d+g+1)(g+1)=9 \leq n+1=\frac{1}{2}H^2+1
\end{equation}
in both cases, proving that $X$ is $BN$ general by Proposition \ref{prop:BN}. 
\end{proof}

Now we can prove

\begin{lemma} \label{lemma:gmeno}
  Let $4 \leq n \leq 9$, $d \geq 1$ and $ g \geq 0$ be integers satisfying 
$g < d^2 /4n+1$.
\begin{itemize}
  \item[(a)] If $g \leq  d^2 /4n$ or if $(n,g) \in \left \{ (5,(d^2+4)/20), (7,(d^2+3)/28) \right\}$, then any Picard minimal $(n,d,g)$-surface is $BN$ general.
\item[(b)] If $g > d^2 /4n$ and $(n,g) \not \in \left \{ (5,(d^2+4)/20), (7,(d^2+3)/28) \right\}$,
then no $(n,d,g)$-surface is $BN$ general.
  \end{itemize}
\end{lemma}

\begin{proof}
  Let $d_0$ and  $g_0$ be as in Lemma \ref{small}. If  we are in case (i) of that lemma, we must have $g \leq  d^2 /4n$ and we 
are done by Corollary \ref{cor:small}. Assume henceforth that we are in case (ii) of Lemma \ref{small}.

If $d_0 \leq n$, set $C_1:=C_0$; if $d_0>n$, set $C_1:=C'=H-C_0$ (as in Lemma \ref{lemma:caso2}). Then, by Lemmas \ref{small} and \ref{lemma:caso2}, we have
$C_1 >0$, with $C_1.H=d_1 \leq n$ and $C_1^2=2(g_1-1) \geq 0$, where
\begin{equation} \label{prodell}
 (d_1,g_1) := \left\{ \begin{array}{ll}
             (d_0,g_0) & \mbox{ if } \; d_0 \leq n, \\
             (2n-d_0,n-d_0+g_0)                 & \mbox{ if } \; d_0>n. \\
             \end{array}
    \right.
\end{equation}
Moreover, $d_1^2-4ng_1=d^2-4ng$. It therefore suffices to prove the statements for $(n,d_1,g_1)$-surfaces, where $d_1 \leq n$.

To prove (a), assume that $X$ is a Picard minimal $(n,d_1,g_1)$-surface and 
$g_1 \leq  d_1^2 /4n$. If $g_1=0$, then $X$ is $BN$ general by  Lemma \ref{rational}. We may therefore assume $g_1>0$. 

We have
\begin{equation} \label{g1}
 g_1 \leq d_1^2/4n \leq n^2/4n=n/4 \leq 9/4, 
\end{equation}
whence $g_1 \leq 2$. If equality holds, the only possibilities are $(n,d_1)=(8,8)$ or $(9,9)$, whence $X$ is $BN$ general by Lemma \ref{lemma:spec}. If $g_1=1$, then \eqref{g1} yields $4n \leq d_1^2$, which implies $d_1 >(n+2)/2$, whence 
$X$ is $BN$ general by Lemma \ref{rational}(ii). This concludes the proof of (a).

To prove (b), assume that $X$ is any  $(n,d_1,g_1)$-surface and 
$g_1 >  d_1^2 /4n>0$. We have 
\[
 (H-C_1)^2= 2(n-d_1+g_1-1) 
= \left\{ \begin{array}{ll}
             2(n-d_0+g_0-1) & \mbox{ if } \; d_0 \leq n, \\
             2(g_0-1)                 & \mbox{ if } \; d_0>n, \\
             \end{array}
    \right.
\]
whence $(H-C_1)^2 \geq 0$ (since we are in case(ii) of Lemma \ref{small}). Since $(H-C_1).H=2n-d_1>0$, it follows that $H-C_1 >0$ by Riemann-Roch. We have
\begin{equation}
  \label{eq:prodII}
  \left(\frac{1}{2}(H-C_1)^2+2\right)\left(\frac{1}{2}C_1^2+2\right) =
(n-d_1+g_1+1)(g_1+1) \geq (g_1+1)^2.
\end{equation}
By Proposition \ref{prop:BN}, in order for $X$ to be $BN$ general, we must have $(g_1+1)^2 \leq \frac{1}{2}H^2+1 =n+1 \leq 10$, whence $g_1 \leq 2$. If $g_1=2$, we must have $n=d_1 \geq 8$ for $X$ to be $BN$ general. But then $d_1^2/4n=n/4 \geq 2=g_1$, a contradiction. Finally, if $g_1=1$, we have $d_1^2 <4n$, which is easily seen to imply that $d_1 \leq (n+2)/2$ for $4 \leq n \leq 9$, except precisely for the two cases $(n,d_1)=(5,4)$ and $(7,5)$. In these two cases we have $d^2-4n(g-1)=d_1^2-4n(g_1-1)=d_1^2$, whence $g=(d^2+4)/20$ and $g=(d^2+3)/28$, respectively. Hence, in all cases in (a), we have $d_1 > (n+2)/2$, and the result follows from Lemma \ref{rational}(ii). Conversely, in all cases in (b) we have 
$d_1 \leq (n+2)/2$, whence 
$X$ is not $BN$ general, again by Lemma \ref{rational}.
\end{proof}

\renewcommand{\proofname}{Proof of Theorem \ref{bncurves}}

\begin{proof}
  Lemmas \ref{primocaso} and \ref{lemma:gmeno} show that for $4 \leq n \leq 9$, the surfaces from Theorem \ref{mainthm} that are $BN$ general are precisely the ones from (i) with $k=1$, from (ii) with $(n,g) \in \left \{ (5,(d^2+4)/20), (7,(d^2+3)/28) \right\}$, and from (iii) and (iv). 
\end{proof}

\renewcommand{\proofname}{Proof}

\section{Proof of Theorem \ref{result}}
\label{sec:proofthm}

We recall the main setting and result  of \cite{kn-TAMS}.

Let $P$ be a smooth projective variety of dimension $r \geq 4$ and
$\E$ a vector bundle of
rank $r-3$ on $P$ that splits as a direct sum of line bundles
\[ \E = \+_{i=1}^{r-3} \M_i. \]
Let
\[ 
s_0=s_{0,1} \+ \cdots \+ s_{0,r-3} \in H^0 (P, \E)= \+_{i=1}^{r-3} H^0(P,\M_i) 
 \]
be a regular section, where $s_{0,i} \in H^0(P,\M_i)$ for $i=1,\ldots,r-3$. 
Set
\[ Y = Z (s_0) \; \; \mbox{and}  \; \; Z=Z(s_{0,1} \+ \cdots \+ s_{0,r-4}) \]
(where $Z=P$ if $r=4$).

Let $X \subset Y$ be a smooth, regular surface (i.e. $H^1(X, \O_X)=0$) and 
$\L$ a line bundle on $X$.

Assume the following: 

\begin{itemize}
\item [(A1)] $Y$ has trivial canonical bundle;
\item [(A2)] $Z$ is smooth along $X$ and the only singularities of $Y$ which lie in $X$ are $\ell$
  nodes $\xi_1, \ldots, \xi_{\ell}$.  Furthermore
  \[ \ell \geq \dim |\L| +2; \]
\item [(A3)] $|\L| \neq \emptyset$ and the general element of $|\L|$ is a smooth, irreducible curve; 
\item [(A4)] for every $\xi_i \in S:=\{\xi_1, \ldots, \xi_{\ell}\}$, if 
$|\L \* \I_{\xi_i}| \neq \emptyset$, then its general member is nonsingular at 
$\xi_i$; 
\item [(A5)] $H^0(C, \N _{C/X}) \iso H^0(C, \N _{C/Y})$ for all
            $C \in |\L|$;
\item [(A6)] $H^1(C, \N _{C/P})=0$ for all $C \in |\L|$;
\item [(A7)] the image of the natural restriction map 
\[  
\xymatrix{
H^0(P, \M_{r-3}) \ar[r] & H^0(S, \M_{r-3} \* \O_S) \iso \CC^{\ell} 
}
\]
has codimension one.
\end{itemize}

Let $s \in H^0(P,\E)$ be a general section.
 
\begin{thm} \label{holger}
 {\rm (\cite[Thm.~1.1]{kn-TAMS})} Under the above setting and assumptions (A1)-(A7), the members of
  $|\L|$ deform to a length $\ell-2 \choose \dim |L|$ scheme of curves
  that are smooth and isolated in the general deformation $Y_t=Z(s_0+ts)$ of $Y_0=Y$. In particular, $Y_t$ contains a smooth, isolated curve that is a deformation of a curve in $|\L|$. 
\end{thm}

We will apply the theorem to prove Theorem \ref{result}.

Let $X$ be a $K3$ surface of degree $2\mu-2$ in $\PP ^{\mu}$ 
that is a complete intersection of type
$(a_1, \ldots ,a_{r-2})$ in some smooth $r$-dimensional ambient space 
$P$, for $r \geq \mu$.
 We will arrange indices so that $a_i \geq 2$ for $i \leq r-4$
and $a_{r-3} \geq a_{r-2}$.  Let
\[ b_i = a_i \; \; \mbox{for} \; \; i=1, \ldots, r-2, \; \; 
\mbox{and} \; \; b_{r-3}=a_{r-3}+a_{r-2}. \]
Then each $b_i \geq 2$ and 
we can construct a Calabi-Yau threefold $Y$ that is a complete intersection 
of type
$(b_1, \ldots, b_{r-3})$ in $P$ as follows: Choose
generators $g_i$ of degrees $a_i$ for the ideal of $X$, so that 
$X = Z(g_1, \ldots, g_{r-2})$. For general 
$\alpha_{ij} \in H^0(P,\O_{P}(b_i-a_j))$ define
\[ f_i:= \sum \alpha_{ij}g_j \]
and
\[ Y := Z(f_1, \hdots, f_{r-3}). \]
If the  
$\alpha_{ij}$'s are chosen 
in a
sufficiently general way, $Y$ has only $\ell= (2\mu-2)a_{r-3}a_{r-2}$ ordinary double points  and they all lie on $X$.  
This can be checked using Bertini's theorem. In fact, the $\ell$ nodes are the intersection points of two general elements of $|\O_X(a_{r-3})|$ and $|\O_X(a_{r-2})|$ (distinct, when $a_{r-3}=a_{r-2}$).  As above, we denote the set of nodes by $S$.

Moreover, for general $\alpha_{ij}$, Bertini's theorem yields that the fourfold
\[ Z := Z(f_1, \hdots, f_{r-4}) \]
is smooth. (Note that $Z=P$ if $r=4$.)

We  are therefore in the setting of Theorem \ref{holger} with
\[ \E :=  \+ _{i=1}^{r-3} \O_{P}(b_i)  \]
and $\M_{r-3}:= \O_{P}(b_{r-3})= \O_{P}(a_{r-3}+a_{r-2})$. We note that $Y$ is nondegenerately embedded in $\PP^m$, where 
\[
m =
\begin{cases} 
\mu & \; \mbox{if $a_{r-2} >1$,} \\
\mu+1 & \; \mbox{if  $a_{r-3} >1$, $a_{r-2} =1$,} \\
\mu+2 & \; \mbox{if $a_{r-3} =a_{r-2}= 1$.} \\
\end{cases} 
\]
The cases where $P=\PP^r$, that is, when $X$ and $Y$ are complete intersections in a projective space, are summarized in \cite[Table~1]{kn-TAMS}.
The following table summarizes all values of
$a_j$, $b_i$, $\ell$, $\mu$, $r$ and $m$ in the cases where $P$ is the intersection of one of the homogeneous spaces mentioned in the introduction and a linear subspace. 

\vspace{0.4cm}
\begin{center}
\begin{tabular}{|c|c|c|c|c|c|c|} \hline
$r$ & $P$ & $(b_i)$ & $(a_j)$ & $\mu$ & $m$ & $\ell$ \\ \hline \hline

$4$  & $G(2, V^5) \cap (1,1)$ & $(3)$  &
$(2,1)$ & $6$ & $7$ & $20$    \\ \hline

$5$ & $G(2, V^5) \cap (1)$ &  $(2,2)$  &
 $(2,1,1)$ &  $6$ & $8$ & $10$    \\ \hline

 $4$ & $\Sigma ^{10} _{12} \cap (1^6)$ &
 $(2)$ &  $(1,1)$  &  $7$ & $9$ & $12$     \\   \hline

 $4$ & $G(V^6, 2) \cap (1^4)$ &   
 $(2)$  & $(1,1)$  &  $8$ & $10$ & $14$     \\   \hline

$4$ & $\Sigma ^6 _{16} \cap (1,1)$ &
 $(2)$ &  $(1,1)$ &  $9$ & $11$ & $16$             \\   \hline

$4$ & $\Sigma ^5 _{18} \cap (1)$ & 
 $(2)$ &  $(1,1)$  &   $10$ & $12$ & $18$     \\   \hline
\end{tabular}
\end{center}
\vspace{0.4cm}

In particular, we note that $a_{r-3}=a_{r-2}=1$ in all cases but the upper one, where $a_{r-3}=2$ and $a_{r-2}=1$. We will henceforth denote this case by $(\dagger)$.

Let now $X$ and $C$ be as in Theorem \ref{bncurves}(II)-(V), with 
$\Pic X \cong \ZZ[H] \+ \ZZ[C]$ and set
$\L:= \O_S(C)$. We will verify conditions (A2)-(A7) and prove the following analogue of \cite[Prop.~7.2]{kn-TAMS}: 

\begin{prop} \label{prop:checkcond}
  Under the contraints given by Theorem \ref{bncurves}(II)-(V), assume that the 
  $\alpha_{ij}$ are general. If 
\begin{equation}
  \label{eq:cddag}
  g \leq 5 \; \; \mbox{or} \; \; 6 \leq g \leq \min\{13,d-6\}  \; \; \mbox{ in case $(\dagger)$}
\end{equation}
and 
\begin{equation}
  \label{eq:cd}
 g \leq \mu-2 \; \; \mbox{otherwise}, 
\end{equation}
the conditions (A1)-(A7) are satisfied.
\end{prop}
 
To prove this, we need some preliminary results. The most delicate condition to check is (A5). By \cite[Rem.~5.1]{kn-TAMS} this condition is equivalent to the condition:
\begin{itemize}
\item [(A5)'] The set of nodes 
$S$ imposes independent conditions on $|\L|$ and the natural map 
$\xymatrix{ \gamma_C: \; \; H^0(C,\N_{X/Y}\* \O_C)  \ar[r] & H^1(C,\N_{C/X})}$
is an isomorphism for all $C \in |\L|$. 
\end{itemize}
For the definition of the map $\gamma_C$ we refer to \cite[(4.4)]{kn-TAMS}.

The following result gives criteria for the first of the two conditions in (A5)'(and (A2))  to hold:

\begin{lemma} \label{lemma:tappato}
Assume that the $\alpha_{ij}$ are general, that $h^0(X,C-H)=0$ and that
$g \leq 13$ in case $(\dagger)$ and $g \leq \mu-2$ otherwise. 
Then $S$ imposes independent conditions on $|\L|$ and condition (A2) holds. 
\end{lemma}

\begin{proof}
  It is easily checked that the proof of \cite[Lemma 6.3]{kn-TAMS} goes through unaffected, which yields that $S$ imposes independent conditions on $|\L|$
if the $\alpha_{ij}$ are general, $h^0(X,\L \* \O_X(-a_{r-2}))=0$ and 
\[
a_{r-2}(2a_{r-3}-a_{r-2})(\mu-1) \geq
\begin{cases} 
g+2 & \; \mbox{if $a_{r-3} \neq a_{r-2}$;} \\
g+1 & \; \mbox{if $a_{r-3}=a_{r-2}$.}
\end{cases} 
\]
Inserting $a_{r-2}
=1$ in all cases, and $a_{r-3}=2$ in case  $(\dagger)$, $a_{r-3}=1$ otherwise, yields the desired numerical conditions. 

Finally, (A2) holds if $\ell \geq g+2$, since $\dim |\L|=g$. This is verified under the above numerical conditions by a direct check. 
\end{proof}

\begin{rem}
  {\rm The last lemma and \cite[Lemma 6.3]{kn-TAMS} only give} sufficient {\rm conditions for $S$ to impose independent conditions on $|\L|$. In \cite{Yu} X.~Yu develops new, interesting methods to check whether $S$ imposes independent conditions on $|\L|$. As an application, he proves the existence of smooth, isolated curves in Calabi-Yau complete intersections threefolds of degrees and genera that are not covered by the existence result \cite[Thm.~1.2]{kn-TAMS}. It is likely that his methods can also be used to improve Theorem \ref{result}.}
\end{rem}

To check the second condition in (A5)', consider the standard short exact sequence

\begin{equation} \label{eq:tangent}
\xymatrix{
0  \ar[r] &  \T _X  \ar[r] & \T _P \* \O_X \ar[r] &  \N _{X/P} \ar[r] & 0
}
\end{equation}

The key observation is the following  result, which is an extension of
\cite[Prop. 1.7]{kl1} to include the cases where the ambient space $P$
is not a projective space:

\begin{prop} \label{key}
  Let $X$ be a $BN$ general $K3$ surface in $\PP^{\mu}$,
  $3 \leq \mu\leq 10$, and $P$ its corresponding ambient space where
  it is a complete intersection.
  The sequence
\[
\xymatrix{     
H^0 (X,\N _{X/P}) \ar[r]^{\hspace{0.3cm}\delta} & H^1 (X,\T_X) \ar[r]^{c} & H^2  (X,\O _X),
}
\]
where $\delta$ is the connecting homomorphism of (\ref{eq:tangent}) and $c$ is the Yoneda
    pairing $ \x c_1(\O_X(1))$, is exact.
\end{prop}

\begin{proof}
Let 
\[
\xymatrix{     
H^0 (X,\N _{X/\PP^{\mu}}) \ar[r]^{\hspace{0.3cm}\overline{\delta}} \ar[r] & H^1 (X,\T_X) 
}
\]  
be the connecting homomorphism of the exact sequence corresponding to \eqref{eq:tangent}
with $P$ replaced by $\PP^{\mu}$. We have $\im \overline{\delta}=\ker c$ by \cite[Prop. 1.7]{kl1}.

The space $H^1 (X,\T_X)$ is the moduli space of all (analytic) infinitesimal first-order deformations of
  $X$, which is well-known to be $20$-dimensional, cf. e.g. \cite{kod} or 
\cite[VIII, Thm.~7.3]{BPHV}.  The space $H^0 (X,\N _{X/P})$ (resp. 
$H^0 (X,\N _{X/\PP^{\mu}})$) 
is the
  moduli space of all infinitesimal first-order deformations of $X$ in $P$ (resp. $\PP^{\mu}$). The map $\delta$ (resp. $\overline{\delta}$) takes a
  deformation class of $X$ to its corresponding deformation class in
  $H^1 (X,\T_X)$ (see e.g. \cite[p.~96]{HM}). We have $\im \delta \sub \im \overline{\delta}$, and the latter is well-known to be $19$-dimensional, cf. e.g. \cite[Thm.~14]{kod}.
Since the primitively polarized $K3$ surfaces that are not $BN$ general
  form a union of countably many proper closed subsets in the moduli space of principally polarized
  $K3$ surfaces, $\im \delta$ is $19$-dimensional and therefore coincides with $\im \overline{\delta}$. Hence $\im \delta= \ker c$.
\end{proof}

For each $C' \in |\L|$ there is an exact sequence
\begin{equation}
  \label{eq:normal}
\xymatrix{     
0 \ar[r] & \N _{C'/X}  \ar[r] & \N _{C'/P} \ar[r] & \N _{X/P} \* \O_{C'} \ar[r] & 0.
}\end{equation}

\begin{lemma} \label{anal}
  Assume that $\O_X(1)$ and $\L$ are independent in $\Pic X$. Then for
  all $C' \in |\L|$ the composition
\[ \xymatrix{     
\phi: H^0 (X,\N _{X/P})  \ar[r] & H^0 (C',\N _{X/P} \* \O_{C'}) \ar[r] & H^1(C', \N _{C'/X})
}
\]
of the restriction morphism with the connecting homomorphism arising
from (\ref{eq:normal}) is surjective. Furthermore, $\ker \phi$ is
independent of $C' \in |L|$.
\end{lemma}

\begin{proof}
  This is proved precisely as in \cite[Lemma 1.10]{kl1}, using Proposition \ref{key} instead of \cite[Prop. 1.7]{kl1}. 
\end{proof}

\begin{prop} \label{thm:gap}
Suppose that the $\alpha_{ij}$ are general and that the line bundles $\L$ and $\O_X(1)$ are independent in 
$\Pic X$.  Then the map $\gamma_C$ in condition (A5)' is an isomorphism for all $C \in |\L|$. 
\end{prop} 

\begin{proof}
  This follows the lines of the proof of \cite[Prop.~6.5]{kn-TAMS} ad verbatim,  replacing $\PP^r$ by $P$ and using Lemma \ref{anal} in the appropriate place at the very end. (Note that $P$ is projectively normal, see e.g. \cite{GW}, so that all sections of $\O_P(a)$, for $a \geq 1$, are induced by homogeneous forms of degree $a$.) 
\end{proof}

\renewcommand{\proofname}{Proof of Proposition \ref{prop:checkcond}}

\begin{proof}
Conditions (A1) and (A3) are obviously satisfied. 
Condition (A4) is satisfied by \cite[Lemma~6.1]{kn-TAMS} and the lines following its proof. 

Condition (A7) is satisfied by \cite[Lemma~6.2]{kn-TAMS}, which goes through unaffected when substituting $\PP^r$ with $P$. 

Condition (A2) is satisfied by Lemma \ref{lemma:tappato}, because of 
\eqref{eq:cddag} and \eqref{eq:cd}. 
 
We next consider condition (A6). Since $X$ is a complete intersection of type $(a_1,...,a_k)$ in $P$,
we have $\N _{X/P} \* \O_{C'} \iso \+ \O _{C'} (a_i)$. Lemma \ref{anal}
together with \eqref{eq:normal} gives
\[ H^1 (C',\N _{C'/P}) \iso H^1 (C',\N _{X/P} \* \O_{C'}) \iso \+ H^1
(C', \O _{C'} (a_i)). \]
In particular, as $\min\{a_i\}=a_{r-2}=1$ in all cases, 
\[ H^1 (\N_{X/P} \* \O_{C'})=0 \hs \mbox{ if and only if } \hs
  H^1 (C', \O _{C'} (1))=0. \]
By \cite[Prop 1.3]{kn} we get that
(A6) is satisfied if and only if
\begin{equation}
  \label{eq:min}
  d \leq 2(\mu-1) \; \; \mbox{or} \; \; d \geq g+\mu 
\end{equation}
holds.  In case $(\dagger)$, this condition reads $d \leq 10$ or $d \geq g+6$, which is satisfied because of \eqref{eq:cddag}. In the other cases, \eqref{eq:cd} implies that if $d > 2(\mu-1)$, then $d >2\mu-2\geq \mu+g$, so that \eqref{eq:min} is again satisfied.

Finally we consider condition (A5), or
equivalently, condition (A5)'.

By Proposition \ref{thm:gap} and the fact that
$\L$ and $\O_X(1)$ are independent in $\Pic X$, the second of the two conditions in (A5)' is satisfied. 
Next we note from the cohomology of
\[ 
\xymatrix{
0  \ar[r] &  \L^{\vee} \ar[r] &  \O_X \ar[r] & \O_C \ar[r] & 0
}
\]
twisted by $H$, Kodaira vanishing and Serre duality, that
\[ h^0(X, C-H) = h^1(\O_C(1)), \]
so that also $h^0(X, C-H)=0$ by conditions \eqref{eq:cddag} and \eqref{eq:cd}, as we have just seen. Thus, by Lemma \ref{lemma:tappato} and the conditions \eqref{eq:cddag} and \eqref{eq:cd}, the first of the two conditions  in (A5)' is satisfied.
Hence, both conditions in (A5)' are satisfied, and so is (A5).
\end{proof}

Theorem \ref{result} is now a consequence of Proposition \ref{prop:checkcond} and Theorem \ref{holger}, recalling that $(n,m,\mu)=(5,7,6)$ in case $(\dagger)$ and $m=\mu+2=n+3$ in the remaining cases. We leave it to the reader to verify that the numerical conditions on $d$ and $g$ in Theorem \ref{result} guarantee that the numerical conditions in both Theorem \ref{bncurves} and in \eqref{eq:cddag} and \eqref{eq:cd} are verified.

\vspace{.5cm}

\end{document}